
\magnification=1200
\hsize=13cm
\vsize=18.5cm
\hoffset=0.2cm
\baselineskip=11pt

\input amssym.def
\input amssym.tex

\font\teneusm=eusm10    
\font\seveneusm=eusm7  
\font\fiveeusm=eusm5    
\newfam\eusmfam

\textfont\eusmfam=\teneusm \scriptfont\eusmfam
=\seveneusm
\scriptscriptfont\eusmfam=\fiveeusm


\def\scr#1{{\scriptstyle{#1}}}

\def\r#1{{\rm #1}}

\def\B#1{{\Bbb #1}}

\footline{{\hss{\vbox to 1cm{\vfil\hbox{\tenrm\folio}}}\hss}}

\centerline{\bf THE TWIN PRIME CONJECTURE}
\bigskip
\centerline{By Yoichi Motohashi}
\bigskip
$$
\matrix{\hbox {The conjecture }\cr
\hbox{`{\it there should be infinitely many
pairs of primes $\{p,p+2\}$}'}\cr
\hbox{ has not been conquered yet. }
}
$$
However, a spectacular drama is now unfolding itself in the theory 
of the distribution of primes.
The complete resolution of the conjecture is thus within the range of modern
mathematics --- {\it perhaps\/}.
Luckily enough, I have been witnessing the series of recent great events 
as a contemporary specialist.
The purpose of the present expository talk is to share 
my excitement with my audience.
Any mathematical discovery is an eventual outcome 
of the rich and long history of our cherished
discipline, and the recent amazing discovery 
by Y. Zhang is a typical instance. 
I shall describe the essence of the fundamental ideas initiated by
GPY (D.A. Goldston, J. Pintz and C.Y. Yildirim) and 
others which had prepared the way for
the discovery, while briefly reviewing the relevant history. 
You will find all basic ideas are 
so simple that you will certainly
be persuaded that the proverb ``{\it small things stir up great\/}" is 
indeed a truth.
\smallskip
Looking back almost half a century ago, I  (then in 
my 20's) was eager to learn Yu.V. Linnik's and
A. Selberg's works in analytic number theory, dreaming
the way to the Never-Never Land of
prime numbers. They taught me
a lot, and I owe them tremendously. I am really happy
that their mathematical spirit is still vividly felt in recent
developments. Indeed, so many wonders
in analytic number theory can be traced back to their ideas.
By trekking further and steadily along the way they prepared, 
you will (I believe) be able to bring us more wonders on primes.
\smallskip
I shall have to be brief in some sections, in order to acquire 
time for more recent work done 
by T. Tao and J. Maynard
independently, which has made Sections 10 and 11 somewhat less
relevant to our main issue
of finding infinitely often bounded differences between primes.
Nevertheless, you will be better off knowing all the facts that
I have put in this text, which I hope will encourage you to
delve into the professional literature on primes.  
\bigskip
\noindent
{\bf Remark 1}: The present text is a substantially improved and augmented
version of the one that had been prepared for my talk delivered at
the Annual Meeting of the Mathematical Society of Japan (15 march 2014).
The expressions that I shall use, whilst being 
adequate for my present (didactic) purpose, are  not
always perfectly precise/correct.
All facts and details on sieve method
and distribution of primes which are needed to
understand recent developments are available in my books [10][12].
\smallskip
\noindent
{\bf Remark 2}: It is highly recommended to visit T. Tao's excellent blog:
\par
\noindent
{\narrower http://terrytao.wordpress.com/2013/06/03/
the-prime-tuples-conjecture\hfill\break-sieve-theory-and-the-work-
of-goldston-pintz-yildirim-motohashi-pintz-and-zhang/\par
}
\par
\noindent
which has various links to more recent developments.
\smallskip
\noindent
{\bf Remark 3.} Because of the digital format specification imposed by arXiv,
two diagrams, one of which was kindly put at my disposal by
the authors of [15], are not included here. 
To view the diagrams,
visit my web-page and download the file EXP2014.pdf.
\smallskip
\noindent
{\it Acknowledgments\/}. I am deeply grateful 
to A. Ivi\'c, M. Jutila, J. Maynard, A. Perelli, T. Tao, N. Watt and H.M. for
their kind comments on the draft of the present text;  
to D.A. Goldston, J. Pintz and C.Y. Yildirim for having been
sharing their epoch-making manuscripts.
\bigskip
\noindent
{\bf 1. The conjecture.}
\par
\noindent
Let
$$
\varpi(n)=\cases{1& $n$ is a prime,\cr 0& $n$ is not a prime,}
$$
and put
$$
\pi(x)=\sum_{n<x}\varpi(n),\qquad
\pi_2(x)=\sum_{n<x}\varpi(n)\varpi(n+2).
$$
Anyone who loves mathematics knows
$$
\pi(x)\sim {x\over \log x}.
$$
Anyone who ardently loves analytic number theory is bitterly defied by the 
conjecture
$$
\pi_2(x)\sim C_0{x\over (\log x)^2},\quad
\hbox{$C_0$: an absolute constant,}\eqno(1.1)
$$
and even by the far more modest statement
$$
\hbox{{\it The twin prime conjecture\/}:\quad
$\displaystyle\lim_{x\to\infty}\pi_2(x)=\infty$.}\eqno(1.2)
$$
\bigskip
\noindent
{\bf 2. To detect twins.}
\par
\noindent
There are two naive means to detect twin primes:
$$
\eqalign{
&(A)\quad \varpi(n)\varpi(n+2)>0,\cr
&(B)\quad \varpi(n)+\varpi(n+2)-1>0.
}
$$
These are of course equivalent to each other as far as
one applies them to {\it individual\/} 
$n$'s, but they are {\it statistically\/} different: always $\varpi(n)
\varpi(n+2)\ge0$ but almost always 
$\varpi(n)+\varpi(n+2)-1=-1$. It appears that opinions 
of sieve specialists are now converging upon
$$
\eqalign{
&\hbox{$(A)$ is too strict,}\cr
&\hbox{$(B)$ is more flexible.}
}
$$
But why?  It is hard to explain the real situation 
to people who are not familiar with sieve method. 
Thus, let me put it bluntly: 
$(A)$ is too exact as it gives 
the definition of $\pi_2(x)$. A sage (M.J.) in analytic number theory 
said that exact formulas contain often too much noise. There were
a lot of attempts, probably since A.M. Legendre's time (the late 18th century),
to clinch to $(1.2)$ by means
of $(A)$; but all eventuated in failure.
In fact, GY (Goldston and Yildirim) commenced their 
investigations in 1999 still brandishing the sharp sword $(A)$. Only in 2004/5,
after a few futile (but highly interesting) attempts with $(A)$, did they 
turn instead to $(B)$. This
was a great turning point in their work. Note that GY actually considered 
{\it primes in tuples\/}:
see Section 7. Here I employ an {\it over\/}-simplification
in order to make the issue clearer. As far as I know, A. Selberg (1950)
was the first who exploited $(B)$, but in a configuration 
different to GY's. 
\bigskip
\noindent
{\bf 3. Sieving out noise.}
\par
\noindent
Imitating the definition of $\pi_2(x)$, one might consider
$$
\sum_{n<x}\big(\varpi(n)+\varpi(n+2)-1\big).\eqno(3.1)
$$
If the sum is positive and large, then the conjecture $(1.2)$ will be resolved.
But this argument is awfully absurd, since obviously $(3.1)$ is 
essentially equal to $2\pi(x)-x$, and one can 
utter only the nonsense
$$
(3.1)\sim\; -x\;.\eqno(3.2)
$$
\par
Nevertheless! Things should look 
pretty different if $(3.1)$ is replaced by
$$
\sum_{n<x}\big(\varpi(n)+\varpi(n+2)-1\big)W(n).\eqno(3.3)
$$
Here $W(n)$'s are {\it non-negative} weights. If one succeeds finding a
nice sequence $\{W(n)\}$ such that $(3.3)$ 
tends to positive infinity, then the conjecture $(1.2)$ will be resolved.
This must be, however, an extremely difficult task, since
such $\{W(n)\}$ should yield a 
considerable dumping of the terms `$1$' and simultaneously should
not affect much the situation of $n$ being a twin prime. That is,
$\{W(n)\}$ is preferably to satisfy
$$
\hbox{$W(n)$ is}\;\cases{\ge0& but very small on average,\cr
\hfil1& when $n$ is a large twin prime.}
$$
\bigskip
\noindent
{\bf 4. Lovely lambda's.}
\par
\noindent
In his work mentioned above,
Selberg employed the $\Lambda^2$-sieve, his great invention
(1947). If translated into our present
situation, it becomes:
$$
\hbox{Consider the quadratic form
$\;\displaystyle\sum_{n<x}\bigg(\sum_{d|n(n+2)}
\lambda(d)\bigg)^2$,}
\atop
\hbox{under the side-condition $\lambda(1)=1$ and $\lambda(d)=0$ for
$d\ge D$,}
$$
where $D$ is a parameter to be fixed optimally eventually,
but initially satisfying only $D<x^{1/2-\varepsilon}$ with an arbitrary small
$\varepsilon>0$.
Expanding the squares out and exchanging the order of summation, we get
the main term and the error term. Selberg diagonalised the main term 
in a highly original way (in fact an application of M\"obius inversion) 
and found an explicit expression 
for optimal $\lambda$'s
that minimises the main term. It is expedient to know that 
these optimal $\lambda$'s satisfy
$$
\lambda(d)\sim\mu(d)\left({\log D/d\over\log D}\right)^2,\eqno(4.1)
$$
with $\mu$ being the M\"obius function, and to note that
$$
\nu(n)>2\quad \Longrightarrow\quad \sum_{d|n}\mu(d)(\log d)^j=0,\;
j\le2.\eqno(4.2)
$$
where $\nu(n)$ is the number of prime factors of $n$ which are
different to each other.
Namely, the choice $(4.1)$ is an {\it approximation\/} to $(4.2)$,
which explains the fact that Selberg's $\lambda$'s yield necessary dumping.
\par
We construct, with these {\it quasi-optimal\/} $\lambda$'s,
$$
W(n)=\left(\sum_{d|n(n+2)}\lambda(d)\right)^2\eqno(4.3)
$$ 
to be used in $(3.3)$. We have, with an appropriate $D$,
$$
\sum_{n<x}W(n)\sim C_1{x\over (\log x)^2},
\eqno(4.4)
$$
and
$$
\sum_{n<x}\big(\varpi(n)+\varpi(n+2)\big)W(n)
\sim C_2 {x\over (\log x)^2}\eqno(4.5)
$$
with certain constants $C_1,\, C_2>0$. Amazing! Compare these with
the conjecture $(1.1)$.
\par
It should be noted that Selberg (ca.\ 1950) examined 
also the use of the weights
$$
\left(\sum_{\scr{d_1|n,\,d_2|(n+2)}\atop
\scr{d_1d_2<D}}\lambda(d_1,d_2)\right)^2,\eqno(4.6)
$$ 
but in a configuration different to $(4.4)$--$(4.5)$ that
 I briefly mentioned already in Section 2.
\vfill\break
\noindent
{\bf 5. RH vs. statistics.}
\par
\noindent
The assertion $(4.5)$ is, in fact, a consequence of
$$
\eqalign{
&\qquad\qquad\hbox{{\it The mean prime number theorem\/}} \cr
&\hbox{For each $A>0$ there exists a $\vartheta>0$ (the {\it level\/})
such that}\cr
&\qquad\sum_{q\le x^\vartheta}\max_{(a,q)=1}\left|\pi(x;a,q)
-{\r{li}(x)\over\varphi(q)}\right|\ll x(\log x)^{-A},\cr
&\qquad\quad\hskip1cm\pi(x;a,q)=\sum_{\scr{n<x}\atop
\scr{n\equiv a\bmod q}}\varpi(n),
}\eqno(5.1)
$$
where $\r{li}$ is the logarithmic integral, $\varphi$ the Euler totient
function, and `$\ll$' means that
the left side is less than a constant multiple of the
right side. The reason why we need
this is simple: With $(4.3)$,
$$
\sum_{n<x}\varpi(n)W(n)
=\sum_{d_1,d_2<D}\lambda(d_1)\lambda(d_2)
\sum_{\scr{a\bmod [d_1,d_2]}\atop
\scr{a(a+2)\equiv0\bmod [d_1,d_2]}}\pi(x;a,[d_1,d_2]),
$$
where $[d_1,d_2]$ is the least common multiple of $d_1,\,d_2$. 
The replacement of each $\pi(x;a,[d_1,d_2])$ by $\r{li}(x)
/\varphi([d_1,d_2])$, providing $(a, d_1d_2)=1$, 
causes an error which can be
estimated with $(5.1)$ if $D<x^{\vartheta/2-\varepsilon}$.
\par
The first result that gave an absolute constant $\vartheta>0$ in $(5.1)$ is 
due to A. R\'enyi (1948).
He exploited the ``large sieve'' of Yu.V. Linnik (1941). This 
{\it statistical equi-distribution\/} 
property of primes among arithmetic progressions
to relatively large moduli must
remind you of the extended Riemann hypothesis ERH. R\'enyi's
prime number theorem states that in some important applications
the extended {\it quasi\/}-Riemann hypothesis could be avoided!
Because of this, a lot of notable people poured their
strenuous efforts into improving upon
R\'enyi's prime number theorem, and E. Bombieri  (1965) established that
$$
\hbox{$(5.1)$ holds for any $\vartheta<{1\over2}$,}\eqno(5.2)
$$
which, {\it in practice\/}, is essentially at the same depth as ERH (actually he 
proved $(5.1)$ with $x^{1/2}(\log x)^{-B(A)}$ in place of $x^\vartheta$). 
I should stress that A.I. Vinogradov (1965) 
proved $(5.2)$ independently; he exploited another 
fundamental innovation due to Linnik: the
``dispersion method'' (1958). 
\par
By the way, in January 1970
I left for Budapest aspiring
to learn analytic number theory under R\'enyi and P. Tur\'an, but
R\'enyi passed away a day after my arrival (1 February).
\bigskip
\noindent
{\bf 6. Powerful modesty.}
\par
\noindent
However, with the best effort one could achieve
only $C_2<C_1$ in $(4.4)$--$(4.5)$. That is,
the asymptotic value $(C_2-C_1)
x/(\log x)^2$ thus attained for $(3.3)$ is 
negative and large, and so is of no more use 
to us than the nonsense $(3.2)$. In fact, in order to truly appreciate
$(4.4)$--$(4.5)$ you ought to be well versed in the theory of the
distribution of primes in arithmetic progressions as well as in sieve method.
Here, be simply amazed that despite its inability to yield
anything about the conjecture $(1.2)$
the assertion comes close to 
the dreamy asymptotic formula $(1.1)$ at least outwardly, 
and moreover, there
we have $W(n)=1$ whenever $n$ is a large twin prime. That is, twin
primes are probably counted in $(4.5)$ but only in an ineffective way;
they must be buried in rubbish.
Then, how to make $(3.3)$ more effective and salvage primes proximate
to each other?  
\par
That is very difficult. The accumulation of past futile attempts suggests
that we ought not to be so daring as to confront $(1.2)$ directly. 
The strategy GY (2004/5) chose was this: We should be modest.
Let us give up trying to directly touch the `twin'. Let us consider instead
$$
\sum_{n<x}\bigg(\sum_{j=1}^k \varpi(n+h_j)-1\bigg)W(n),\eqno(6.1)
$$
with a new $\{W(n)\}$.
Here $h_1<h_2<\cdots<h_k$ are even integers. They should not be
trivial like $\{2,4,6\}$ because one of $n+2,n+4,n+6$ is always
divisible by $3$. A natural condition
on the tuple $\{h_1,h_2,\ldots,h_k\}$ is that
$$
\hbox{the number of different $h_j\bmod p$ be less than $p$ for any prime
$p$,}\eqno(6.2)
$$
which avoids the redundancy that a member among $\{n+h_j: j=1,\ldots, k\}$
is always divisible by a fixed prime. Obviously,
$$
\hbox{$\displaystyle\sum_{j=1}^k \varpi(n+h_j)-1>0$}
\atop\hbox{$\Longrightarrow\;\{n+h_1,n+h_2,\ldots, n+h_k\}$ 
contains at least two primes.}
$$
If this holds with infinitely many $n$, then
$$
\liminf_{t\to\infty}\,(p_{t+1}-p_t)\le h_k-h_1,
$$
with $p_t$ the $t$-th prime. Bounded differences between primes should occur 
infinitely often. The establishment of this will be a tremendous achievement,
even though it is perhaps less impressive than the ultimate assertion $(1.2)$. 
\bigskip
\noindent
{\bf 7. Gem box principle.}
\par
\noindent
We have to choose the weights $\{W(n)\}$ in $(6.1)$.
Here a truly decisive observation was made by
GPY (2005): Let
$P(n)=(n+h_1)(n+h_2)\cdots(n+h_k)$. Then,
$$
\hbox{$\nu(P(n))=k+\ell$ with $
0\le \ell<k$}\quad\Longrightarrow\quad
{\hbox{there are at least $k-\ell$ primes}\atop
\hbox{among $n+h_1,\ldots,n+h_k$}\,.} \eqno(7.1)
$$
This is an application of Dirichlet's pigeon box principle; but I 
very much prefer {\it gems\/} to pigeons. Here  
$n$'s are actually to be restricted so that $(7.1)$ is valid, 
which can be realised in a simple way that does not cause
any loss of generality. 
\bigskip
\noindent
{\bf 8. Magical tapering.}
\par
\noindent 
The new parameter $\ell\ge 0$ is
to be incorporated. In practice, however, it is hard to utilise
$(7.1)$ without making any compromise; that would be
a return to the stiffness we wished to depart from. I am not very sure if this is
what really occurred to them, but GPY seem to have
turned to Selberg's argument
which I indicated in the first paragraph of Section 4. The relevant
approach is to consider
$$
\sum_{n<x}\bigg(\sum_{d|P(n)}\lambda(d)\bigg)^2,
\quad\cases{\lambda(1)=1,\cr\lambda(d)=0,&$d>D.$}
$$
The optimal $\lambda$ satisfies
$$
\lambda(d)\sim
\mu(d)\left({\log{D/d}\over\log D}\right)^k.
$$
Then, GPY practised real magic by introducing
$$
\hbox{the further tapering factor 
$\displaystyle\left({\log D/d\over\log D}\right)^\ell$},
$$
and they constructed the weight
$$
W(n)=\left(\sum_{d|P(n),\,d<D}
\mu(d)\left({\log D/d\over\log D}\right)^{k+\ell}\right)^2.\eqno(8.1)
$$
As a matter of fact,
this is an approximation to the filtering concerning $(7.1)$, since
$$
\nu(P(n))> k+\ell\quad
\Longrightarrow\quad\sum_{d|P(n)}
\mu(d)(\log d )^j=0,\; j\le k+\ell.
$$
\bigskip
\noindent
{\bf 9. Divine multiplier.}
\par
\noindent
With $W(n)$ as in $(8.1)$, GPY  
computed asymptotically the sums
$$
\eqalign{
T^{(1)}(x;k,\ell;D)&=\sum_{n<x}W(n),\cr
T^{(2)}(x;k,\ell;D)&=\sum_{n<x}
\bigg(\sum_{j=1}^k\varpi(n+h_j)\bigg)W(n).
}\eqno(9.1)
$$
They discovered that, with $D=x^{\vartheta/2}$ 
($\vartheta$ as in $(5.1)$) and a positive $\Delta(x)
\approx x(\log x)^{-k}$, one has:
$$
\eqalignno{
\Big(T_P^{(2)}-T_P^{(1)}\Big)(x;k,\ell;D)&\sim
\Big(\vartheta\cdot{k\over k+2\ell+1}\cdot
{2\ell+1\over \ell+1}-1\Big)\Delta(x).&(9.2)
}
$$
This multiplier of $\Delta(x)$ is probably {\it
one of the greatest surprises\/} in the entire history of
number theory. Setting $\ell=[\sqrt{k}]$ for instance, we find readily that
$$
\hbox{if $\vartheta>{1\over2}$ and $k$ large, then
$\{n+h_1,n+h_2,\ldots,n+h_k\}$}\atop
\hbox{ contains at least two primes. $\Longrightarrow$
Bounded differences between primes!}\eqno(9.3)
$$
If you had not the extra parameter $\ell$; that is, if you put $\ell=0$,
then $(9.2)$ would be nothing. Without $\vartheta>1$, which is 
truly beyond any science fiction, 
nothing would come out from $(9.2)$ with $\ell=0$. In fact it is
known that $(5.1)$ does not hold for any $\vartheta>1$.
\bigskip
\noindent
{\bf 10. Divide and conquer.}
\par
\noindent
The assertion $(9.3)$ is indeed wonderful, if only
one can leap beyond the barrier $\vartheta={1\over2}$ in 
$(5.1)$.
\smallskip
Let me be a little bit personal: I may count myself as one of the earliest
people who tried seriously to make this leap, of course without any surmise of 
recent developments. I was aware at least that
not the large sieve but the dispersion method of Linnik is the key. But
I could publish only a short report (1976) which relied still on the
large sieve; my work relevant to the dispersion method was 
utterly incomplete, which was inevitable because of my meagre experience
with the theory of exponential sums \`a la A. Weil.
 Later BFI (Bombieri, 
J. B. Friedlander and H. Iwaniec (1986)) made a remarkable
progress in this direction. Their main 
result is valid with any 
$\vartheta<{4\over7}$, 
but under a restriction on the moduli of the
arithmetic progressions which makes it inadequate for 
the computation of the second sum in $(9.1)$. 
\smallskip
Thus a genuinely new insight was needed 
into the problem $(6.1)$ and the barrier problem. In this situation
an idea occurred to MP (2005) (see [11][14]
as well); actually we each
independently had essentially the same idea, which involved the use of
some corner-cutting in order to break the stalemate. On my side:
soon after getting the
first version of GPY (from G in early April 2005) I realised that a 
 {\it smoothing\/} could be applied to the summation variable
$d$ in $(8.1)$. That is, we need not sum over all $d<D$ but it
suffices to restrict ourselves to those $d$ which have relatively 
{\it small\/} prime divisors only; 
I mean that even after applying such a corner-cutting the multiplier
of $\Delta(x)$ in $(9.2)$ does not change essentially, 
although $\Delta(x)$ itself ought to be altered accordingly. 
\par
Actually, MP (2005/6) modified the argument of 
GGPY (GPY and S. Graham (2005)) in order to
incorporate this smoothing. Let me nevertheless employ an
asymptotic expression
for the sake of temporary convenience.
Then, what MP did is the same as to replace $(8.1)$ by
$$
W(n)=\left(
\mathop{\quad{\sum}^{(\omega)}}
_{d|P(n),\,d<D}
\mu(d)\left({\log D/d\over\log D}\right)^{k+\ell}\right)^2,\eqno(10.1)
$$
where $\sum^{(\omega)}$ indicates that all prime divisors of $d$ are
less than $D^\omega$.
Then the multiplier in $(9.2)$, of course under the new setting, 
is found to be larger than
$$
\vartheta_\r{MP}\cdot{k\over k+2\ell+1}\cdot
{2\ell+1\over \ell+1}-1-\exp(-3k\omega/8),\eqno(10.2)
$$
provided that one has, for any given $A>0$,
$$
\mathop{\quad{\sum}^{(\omega)}}_{q\le x^{\vartheta_\r{MP}}}
\sum_{\scr{(a,q)=1}\atop\scr{P(a)\equiv 0\bmod q}}\left|\pi(x;a,q)
-{\r{li}(x)\over\varphi(q)}\right|\ll x(\log x)^{-A},\eqno(10.3)
$$
where 
$\sum^{(\omega)}$ means that all prime factors of $q$ are
less than $x^\omega$. 
Here I am not very precise, since MP
tacitly assumed for the sake of convenience
that $\ell\approx\sqrt{k}$,
$ \omega\approx 1/\sqrt{k}$ with $k$ large;
however, these assumptions are not of critical importance for the application in
question, that is, to detect infinitely often bounded differences 
between primes. 
I remark also that the hypothetical mean prime number
theorem which is required by MP is a consequence of $(10.3)$;
that is, MP assumed in fact somewhat less. Anyway we have:
$$
\hbox{$\vartheta_\r{MP}>{1\over2}$ in $(10.3)$}\atop
\hbox{ $\Longrightarrow\quad$ bounded differences 
between primes occur infinitely often.}\eqno(10.4)
$$
Why is this important? Because,
with $(10.3)$, instead of $(5.1)$, the feasibility of a proof by
the dispersion method of Linnik becomes much higher. More precisely, 
the smoothing yields a {\it quasi\/}-infinitely 
factorable structure in the moduli set
$\{q\}$;  namely, we now have instead
$$
\{q_1q_2: q_1\le Q_1, q_2\le Q_2\}, 
$$
essentially for any
multiplicative decomposition $Q_1Q_2\le x^{\vartheta_\r{MP}}$. In practice, 
we put the summation over $q_1$, say, 
outside and consider the {\it dispersion\/} 
of the inner sum over $q_2$, via the Cauchy inequality. We will be able to
detect more cancellation than with the ordinary setting $(5.1)$.
Further, we may appeal to R.C. Vaughan's reduction argument (1980),
or the like, in dealing with the sums over primes. 
This strategy is nothing other than ``divide et impera''.
\bigskip
\noindent
{\bf 11. From nowhere.}
\par
\noindent
 As to the proof of $(10.3)$ for a $\vartheta_\r{MP}>{1\over2}$, 
 I was somehow inclined to be optimistic; 
 and I thought I would have `time'. Thus,
in the mean time, I was playing with
automorphic $L$-functions, enjoying some success, 
but for too long perhaps.
Then, in early April last year I felt a jolt. 
The epicentre was an
unknown mathematician named Y. Zhang; I mean that the man had not been
known among specialists. Soon I got a copy of his paper (probably a draft). 
I felt as if I had seen it some 7 years ago, for its overall strategy was the
same as that of MP(2005/6). 
\par
Of course I was truly
impressed by the extremely important
fact that Zhang cleared away the level barrier in the
context of $(10.3)$.
The man who came from nowhere struck the target indeed. 
Therefore, mankind has now
$$
\liminf_{t\to\infty}(p_{t+1}-p_t)< \infty.\eqno(11.1)
$$
\smallskip
To achieve $(10.3)$, for some $\vartheta_\r{MP}>{1\over2}$, 
Zhang appealed to
P. Deligne's famous work (1980) on the Weil conjecture; 
in this respect, he followed, to a large
extent, the work by BFI mentioned above. Thus I am
unable to confirm his reasoning on my own but have to rely on
the affirmative opinion of experts. 
I have no courage to exploit any result which
I do not fully understand; neither have I any other way than to trust, with
considerable caution,
competent authors whose claims depend on works which are far beyond
my expertise. Nevertheless, here I may try to explain why 
such heavy machinery comes into play in dealing with $(10.3)$.
In essence, it is because of the factoring 
of various terms and summation intervals, which is described 
in the previous section. I mean that the strategy there
reduces the problem into pieces, all of which are more or less equivalent to
counting integers in various arithmetic progressions.
To manage this entangled task, presently
we have essentially only one means: the Poisson
summation formula. Main terms are not troublesome, though 
often complicated. Real trouble comes naturally from the tail parts, which
are expressed in terms of finite or infinite exponential sums. 
Arguments of the exponentiated terms involve rational
numbers with varying numerators and denominators;
then Deligne's work becomes relevant, as it gives strong  and {\it uniform\/}
bounds for such sums.
\bigskip
\noindent
{\bf 12. Phase transition.}
\par
\noindent
Another sensation came more recently from a 
postdoc: J. Maynard (November 2013),
claiming
$$
\liminf_{t\to\infty} (p_{t+1}-p_t)\le 600.\eqno(12.1)
$$
What is really sensational is in his statement that
his argument {\it does not incorporate any of the technology used by Zhang; 
the proof is essentially elementary, relying only on the
Bombieri--Vinogradov theorem\/}, i.e., $(5.2)$. This is a true
phase transition, and a great gift
to all who feel uneasiness when they have to chew
works that depend on the highly demanding work of Deligne and 
A. Weil (1949),  even though the efforts of 
S.A. Stepanov (since 1969) have yielded accessible 
elementary proofs of some of the consequences of their work.
\par
And more. According to Maynard, Tao (October 2013) 
got essentially the same idea; and they independently established, only on
R\'enyi's $(5.1)$,
$$
{\hbox{For each $m\ge2$ there exists a $k$ such that}\atop
\hbox{with any $\{h_j\}$ satisfying $(6.2)$}}\atop
{\hbox{the tuple $\{n+h_1,n+h_2,\ldots,n+h_k\}$}
\atop\hbox{ contains at least $m$ primes for infinitely many
$n$.}}\eqno(12.2)
$$
They even got an estimate for $k$
in terms of $m$. Fantastic! 
\par
Their argument is, 
to some extent, a realisation as well as an extension 
of Selberg's approach $(4.6)$. Hence,
in a sense, $(12.1)$ would have been possible to attain
in 1965 when $(5.2)$ was established; and $(12.2)$  in 1950! 
By this I mean that for more than half a century,  indeed
until a few months ago, no sieve experts had ever tried to seriously 
look into the
ending remark (on p.245) in Selberg's `Lectures on sieves'. 
I should of course add that the phase transition brought about 
by Maynard and Tao
was an outcome of the sieve movement commenced 
by Goldston and Yildirim in 1999, without which I suspect that not only 
Maynard--Tao's discovery but also the recent
wonders concerning bounded differences between
primes would have remained under sand, and perhaps
would have lain undiscovered for decades to come.
Better ideas always survive; what I described in the last two sections may
appear obsolete, at least for now.
\par
The key points of Maynard's argument are as follows:
Basically we are dealing with the quadratic form
$$
\sum_{n<x}
\left(\sum_{d_j|(n+h_j),\;\forall j\le k}
\lambda(d_1,d_2,\ldots,d_k)\right)^2,\quad
d_1d_2\cdots d_k\le D.\eqno(12.3)
$$
We need to be cautious in dealing with the prime factors of $d_j$; but let
us ignore this presently: a correct procedure 
is indicated in Appendix below.
Then, in a fashion familiar 
to those who are experienced in dealing 
with sums of arithmetical functions in sieve method, an application 
of Selberg's change of variables (in fact, an instance 
of the M\"{o}bius inversion) allows one to express
$\lambda$'s in terms of any given $F(\xi_1,\xi_2,\ldots ,\xi_k)$ 
as far as $F$ is supported 
on $\{\xi_1+\xi_2+\cdots +\xi_k\leq 1\,:\,\xi_j\geq 0,\forall j\leq k\}$.
This is in fact an extension of the argument due to GGPY (2005); their choice 
corresponds to the specialisation $F\big(\,\underline{\xi}\,\big)=
f(\xi_1+\xi_2+\cdots+\xi_k)$. We let $W(n)$ stand for the
squares in $(12.3)$ with such $\lambda$'s, and
engage in the evaluation of 
$$
\sum_{n<x}\bigg(\sum_{j=1}^k \varpi(n+h_j)-\rho\bigg)W(n),
\eqno(12.4)
$$
which is an obvious analogue of $(6.1)$; the parameter $\rho$ is to be
fixed later. Actually we need to apply
{\it pre-sifting\/} to $n$'s as indicated in $(A.3)$ below, which
is not of absolute necessity but for the sake of technical comfort in
dealing with $d$'s coming from $(12.3)$.
In this way, with $\vartheta$ as in (5.1), 
we find that the appropriate analogue of 
the multiplier of $\Delta(x)$ in (9.2) is:
$$
{\vartheta\over2}\sum_{j=1}^k{J_k^{(j)}}(F)-\rho I_k(F),
\eqno(12.5)
$$
where
$$
\eqalign{
I_k(F)&\,=\int_0^1\!\cdots\!\int_0^1
F(\xi_1,\xi_2,\ldots,\xi_k)^2d\xi_1\cdots d\xi_k,\cr
J_k^{(j)}(F)&\,=\int_0^1\!\cdots\!\int_0^1\left(
\int_0^1F(\xi_1,\xi_2,\ldots,\xi_k)d\xi_j\right)^2d\xi_1\cdots d\xi_{j-1}
d\xi_{j+1}\cdots d\xi_k.
}
$$
If we put $\rho=1$ and
$F\big(\,\underline{\xi}\,\big)=(1-\xi_1-\cdots-\xi_k)^\ell$,
then we recover $(9.2)$ due to GPY (2005).
\par
We are naturally interested in the variation problem
$$
M_k=\sup_F{\sum_{j=1}^k{J_k^{(j)}}(F)\over I_k(F)},
$$
where the supremum is over functions $F$ that
 are piece-wise differentiable in the domain indicated above 
and such that
$I_k(F)\ne0$, $J_k^{(j)}(F)\ne0$ for each $j\le k$. Let
$$
\rho=m-1,\quad m=\inf\{r\in\B{N}: r\ge\vartheta M_k/2\}.
$$ 
Then one finds that there are at least $m$ primes
in $\{n+h_1,n+h_2,\ldots, n+h_k\}$ for infinitely many $n$'s.
With a delicate optimisation, Maynard has found
$$
M_{105}>4.002,
$$
which together with $(5.2)$ implies $(12.1)$ as there 
exists $\{h_1,h_2,\ldots,h_{105}\}$ such
that $h_{105}-h_1=600$. More strikingly, he has shown via a 
simple choice of $F$ that for sufficiently large $k$
$$
M_k>\log k-2\log\log k-2.
$$
This implies $(12.2)$. 
\bigskip
I repeat: R\'enyi established his prime number 
theorem $(5.1)$ in 1948 and the argument of Manynard and Tao
has its root in Selberg's work of 1950. Thus, more than 60 years ago
when I entered elementary school, the notion
that bounded differences between primes occur infinitely often
could easily have already belonged to common knowledge.
\bigskip
\noindent
{\bf Appendix.} As an induction for students who intend 
to study Maynard's work,
I shall provide details of his {\it arithmetic\/}
manipulations in the case $k=2$, which is enough typical so that
one may readily infer that the general case 
is to be settled as is shown in $(12.5)$. As to Tao's
argument, the difference is only 
in the way of computing asymptotically the main terms
which arise after sieving. He employed Fourier analysis in place of 
the usual method of summing arithmetic functions which Maynard used;
see Tao's polymath8 blog, the address of which is
given in the references below.
\smallskip
We assume that $N$ tends to infinity, and we put
$$
Y=\log\log N, \quad Z=\prod_{p\le Y}p.\eqno(A.1)
$$
The r\^ole of $Y$ or rather that of $Z$
 is important, as it makes the co-primality requirement 
 in various sums easy to attain and also yields crucial truncations
after the change of variables in the mode of Selberg; for the latter,
see $(A.10)$, for instance. The prime number theorem implies 
$Z\ll (\log N)^2$, which can be regarded to 
be negligibly small in our discussion.
We choose $c_0\bmod Z$ to satisfy $(Z, (c_0+h_1)(c_0+h_2))=1$,
which is possible whenever $\{h_1,h_2\}$ satisfies 
the case $k=2$ of $(6.2)$.
We shall work on the  assumption:
$$
\hbox{$\lambda(u,v)=0$ if any of the following holds}
\atop
\hbox{$uv>D, \;|\mu(uv)|=0, \; \left(uv,\,Z\right)>1$.}
\eqno(A.2)
$$
With this, we shall consider
$$
\sum_{\scr{N\le n<2N}\atop\scr{n\equiv c_0\bmod Z}}
\left(\sum_{d_1|(n+h_1),d_2|(n+h_2)}
\lambda(d_1,d_2)\right)^2.\eqno(A.3)
$$
Because of the choice of $c_0$ and since $N$ is large, 
we have always $(n+h_1, n+h_2)=1$
and thus $(d_1,d_2)=1$ in $(A.3)$, conforming
with $(A.2)$. We shall exploit this fact in the sequel without mention.
\par
Expanding the squares and 
changing the order of summation, we see that the sum equals
$$
(N/Z)S_0+O\big(\lambda^2_\r{max}(D\log D)^2\big),\eqno(A.4)
$$
where $\lambda_\r{max}=\sup|\lambda(u,v)|$ and
$$
S_0=\sum_{\scr{d_1,f_1, d_2, f_2}\atop
\scr{(d_1f_1,d_2f_2)=1}}{\lambda(d_1,d_2)\lambda(f_1,f_2)\over
[d_1,f_1][d_2,f_2]}.\eqno(A.5)
$$
Because of $(A.2)$, the condition $(d_1f_1,d_2f_2)=1$ is equivalent
to $(d_1, f_2)(d_2,f_1)=1$. Then we have
$$
\eqalignno{
S_0=&\sum_{\scr{d_1,f_1, d_2, f_2}\atop
\scr{(d_1, f_2)(d_2,f_1)=1}}{\lambda(d_1,d_2)\lambda(f_1,f_2)\over
d_1d_2f_1f_2}\sum_{u_1|(d_1,f_1),\,u_2|(d_2,f_2)}
\varphi(u_1)\varphi(u_2)\cr
=&\sum_{u_1,u_2}
\varphi(u_1)\varphi(u_2)
\sum_{{\scr{d_1, f_1, d_2,f_2}\atop
\scr{(d_1, f_2)(d_2,f_1)=1}}\atop\scr{u_1|(d_1,f_1), u_2|(d_2,f_2)}}
{\lambda(d_1, d_2)\lambda(f_1, f_2)\over d_1d_2f_1f_2}\cr
=&\sum_{u_1,u_2}
\varphi(u_1)\varphi(u_2)
\sum_{\scr{d_1, f_1, d_2,f_2}\atop
\scr{u_1|(d_1,f_1), u_2|(d_2,f_2)}}
{\lambda(d_1, d_2)\lambda(f_1, f_2)\over d_1d_2f_1f_2}
\sum_{v_1|(d_1,f_2),\,v_2|(d_2,f_1)}\mu(v_1)\mu(v_2)\cr
=&\sum_{u_1,u_2, v_1,v_2}\varphi(u_1)\varphi(u_2)\mu(v_1)
\mu(v_2)\sum_{{\scr{d_1, f_1, d_2,f_2}\atop
\scr{u_1|(d_1,f_1), u_2|(d_2,f_2)}}\atop\scr{v_1|(d_1,f_2),\, 
v_2|(d_2, f_1)}}
{\lambda(d_1, d_2)\lambda(f_1, f_2)\over d_1d_2f_1f_2}\cr
=&\sum_{u_1,u_2, v_1,v_2}\varphi(u_1)\varphi(u_2)\mu(v_1)
\mu(v_2)\sum_{\scr{u_1v_1|d_1}\atop\scr{u_2v_2|d_2}}
{\lambda(d_1,d_2)\over d_1d_2}
\sum_{\scr{u_1v_2|f_1}\atop\scr{u_2v_1|f_2}}
{\lambda(f_1,f_2)\over f_1f_2}.&(A.6)
}
$$
Hence, we put
$$
\eta(w_1,w_2)=\mu(w_1)\mu(w_2)
\varphi(w_1)\varphi(w_2)
\sum_{\scr{d_1, d_2}\atop\scr{w_1|d_1,\,w_2|d_2}}
{\lambda(d_1,d_2)\over d_1d_2},\eqno(A.7)
$$
and have
$$
S_0=\sum_{\scr{u_1,u_2,v_1,v_2}\atop
\scr{(u_1u_2v_1v_2,Z)=1}}\mu^2(u_1u_2v_1v_2)
{\eta(u_1v_1,u_2v_2)\eta(u_1v_2,u_2v_1)\over
\varphi(u_1)\varphi(u_2)}\cdot{\mu(v_1)\mu(v_2)
\over(\varphi(v_1)\varphi(v_2))^2}\,.\eqno(A.8)
$$
Applying the M\"obius inversion formula to $(A.7)$, we have
$$
\lambda(d_1,d_2)=\mu(d_1)\mu(d_2)d_1d_2
\sum_{{\scr{w_1,w_2}\atop\scr{(w_1w_2,Z)=1}}\atop
\scr{d_1|w_1,d_2|w_2}}\mu^2(w_1w_2)
{\eta(w_1,w_2)\over\varphi(w_1)\varphi(w_2)}.
\eqno(A.9)
$$
The condition $(A.2)$ is readily seen to be well satisfied
with any {\it any\/} $\eta(u,v)$ as far as it vanishes for $uv>D$. 
Namely, under this specification of $\eta$
one may regard $(A.9)$ as the definition of $\lambda$'s,
as we shall do in the sequel.
Then, $(A.8)$ implies that
$$
S_0=\sum_{\scr{u_1,u_2}\atop\scr{(u_1u_2,Z)=1}}\mu^2(u_1u_2)
{\eta^2(u_1,u_2)\over\varphi(u_1)\varphi(u_2)}
+O\big(\eta_\r{max}^2(\log D)^2/Y\big),\eqno(A.10)
$$
since we have
$$
\sum_{u\le D}{1\over\varphi(u)}\ll\log D,\quad
\sum_{v>1,\,(v,Z)=1}
{1\over\varphi(v)^2}\ll Y^{-1}.
\eqno(A.11)
$$
\par
Next, we shall consider
$$
\sum_{\scr{N\le n<2N}\atop\scr{n\equiv c_0\bmod Z}
}\varpi(n+h_1)\left(\sum_{d_1|(n+h_1),d_2|(n+h_2)
}\lambda(d_1,d_2)\right)^2.\eqno(A.12)
$$
It makes no difference if the condition 
$d_1\mid (n+h_1)$ is replaced by 
the apparently stronger condition 
$d_1=1$, and so we see that $(A.12)$ equals
$$
{1\over\varphi(Z)}
(\r{li}(2N)-\r{li}(N))S_1+O\big(\lambda^2_\r{max}
 E_3(2N,D^2Z)\big),\eqno(A.13)
$$
where
$$
S_1=\sum_{d,f}{\lambda(1,d)\lambda(1,f)\over\varphi([d,f])}
\eqno(A.14)
$$
and
$$
E_l(x,Q)=\sum_{q\le Q}\tau_l(q)
\max_{(a,q)=1}\left|\pi(x;a,q)-{\r{li}(x)
\over\varphi(q)}\right|.\eqno(A.15)
$$
Here $\tau_l(q)$ is the number of ways expressing $q$ as a product
of $l$ factors; in fact, the number of representations of
$q$ as the least common multiple of two integers is bounded by $\tau_3(q)$.
Using th relation
$$
{\varphi(d)\varphi(f)\over\varphi([d,f])}=\sum_{u|(d,f)}\gamma(u),
\quad \gamma(u)=\prod_{p|u}(p-2)
\eqno(A.16)
$$
we have
$$
S_1
=\sum_u\gamma(u)\Bigg(\sum_{u|d}{\lambda(1,d)\over\varphi(d)}
\Bigg)^2.\eqno(A.17)
$$
Imitating $(A.7)$, we put
$$
\eta_1(u)=\mu(u)\gamma(u)
\sum_{u|d}{\lambda(1,d)\over\varphi(d)},\eqno(A.18)
$$
so that
$$
S_1=\sum_u{\eta_1^2(u)\over\gamma(u)}.\eqno(A.19)
$$
Inserting $(A.9)$ into $(A.18)$, we have, after an arrangement,
$$
\eqalign{
\eta_1(u)&=u\gamma(u)\mu(u)\sum_{\scr{(w_1w_2,Z)=1}
\atop\scr{u|w_2}}\mu^2(w_1w_2)
{\eta(w_1,w_2)\mu(w_2)\over\varphi(w_1)\varphi^2(w_2)}\cr
&={u\gamma(u)\over\varphi^2(u)}
\sum_{(w_1u,Z)=1}\mu^2(w_1u)
{\eta(w_1,u)\over\varphi(w_1)}
+O\big(\eta_\r{max}(\log D)/Y\big).
}\eqno(A.20)
$$
This error term is due to the fact that if $w_2\ne u$, then $w_2/u>Y$. 
Further, we have
$$
\eta_1(u)=\sum_{(w_1u,Z)=1}\mu^2(w_1u)
{\eta(w_1,u)\over\varphi(w_1)}+O\big(\eta_\r{max}(\log D)/Y\big),
\eqno(A.21)
$$
since
$$
{u\gamma(u)\over\varphi^2(u)}=\prod_{p|u}\left(
1-{1\over (p-1)^2}\right)=1+O(1/Y),\quad u>1.\eqno(A.22)
$$
\par
With this, we put
$$
\eta(d_1,d_2)=F\left({\log d_1\over \log D},{\log d_2\over\log D}\right),
\eqno(A.23)
$$
where $F$ is as in the last section but with $k=2$. 
Collecting $(A.10)$, $(A.19)$ and $(A.21)$,
we find that we need to evaluate asymptotically the sums
$$
\eqalign{
&\sum_{\scr{u_1,u_2}\atop\scr{(u_1u_2,Z)=1}}
{\mu^2(u_1u_2)\over\varphi(u_1)\varphi(u_2)}
F\left({\log u_1\over \log D},{\log u_2\over\log D}\right)^2,\cr
&\sum_{\scr{u}\atop\scr{(u,Z)=1}}{1\over\gamma(u)}
\left(\sum_{\scr{w_1}\atop\scr{(w_1,Z)=1}}
{\mu^2(w_1u)\over\varphi(w_1)}
F\left({\log w_1\over \log D},{\log u\over\log D}\right)\right)^2.
}\eqno(A.24)
$$
Here one may replace $\mu^2(u_1u_2)$ by $\mu^2(u_1)\mu^2(u_2)$
and do the same with the factor $\mu^2(w_1u)$, since 
$\mu(u_1u_2)=0$, for instance, implies that $u_1$ and $u_2$ 
are divisible by a $u>Y$ and such terms can be discarded in much
the same way as is done in $(A.10)$. Thus,
the computation can be performed 
in a fashion quite familiar in the theory of sums
of arithmetic functions weighted with smooth functions; in essence it is
an application of summation/integration by parts. We may
skip the details and show only the end result: The last two sums are
asymptotically equal to
$$
\eqalign{
& (\log D)^2\left({\varphi(Z)\over Z}
\right)^2\int_0^1\int_0^1F^2(\xi_1,\xi_2)d\xi_1d\xi_2,\cr
& (\log D)^3\left({\varphi(Z)\over Z}
\right)^3\int_0^1\left(\int_0^1F(\xi_1,\xi_2)d\xi_1\right)^2d\xi_2,
}\eqno(A.25)
$$
respectively, as $D$ tends to infinity.
\par
Now, we choose $D=N^{\vartheta/2-\varepsilon}$, 
with $\vartheta$ as in $(5.1)$.
Then, the assertions $(A.4)$ and $(A.13)$ yield the multiplier
$$
\eqalign{
{\vartheta\over2}\Bigg[\int_0^1\bigg(\int_0^1&
F^2(\xi_1,\xi_2)d\xi_1\bigg)^2d\xi_2+
\int_0^1\bigg(\int_0^1F(\xi_1,\xi_2)d\xi_2\bigg)^2d\xi_1\Bigg]\cr
-\rho&\int_0^1\int_0^1F^2(\xi_1,\xi_2)d\xi_1d\xi_2
}\eqno(A.26)
$$
for the sum
$$
\sum_{\scr{N\le n<2N}\atop\scr{n\equiv c_0\bmod Z}}
\big(\varpi(n+h_1)+\varpi(n+h_2)-\rho\big)W(n),\eqno(A.27)
$$
where $W(n)$'s stand for the squares in $(A.3)$ with $\lambda$'s as
in $(A.9)$ along with $(A.23)$. 
We may skip the estimation of the error terms coming from
$(A.10)$ and $(A.21)$ as they should not cause any difficulty. 
As to the error term in $(A.13)$, we need to eliminate the factor
$\tau_3(q)$ in $(A.15)$. This can be achieved via an application of
the Cauchy inequalty; that is,
$$
E^2_l(x,Q)\ll x(\log Q)^{l^2}E_1(x,Q).\eqno(A.28)
$$
\bigskip
\centerline{\bf References}
\bigskip
\item{[1]} E. Bombieri: {\it Le Grand Crible dan la 
Th\'eorie Analytique des Nombres\/} (second \'ed.). 
Ast\'erisque {\bf 18}, Paris 1987. 
\item{[2]} E. Bombieri, J.B. Friedlander and H. Iwaniec: Primes in
arithmetic progressions to large moduli.\ I. Acta Math., {\bf 156}(1986),
203--251.
\item{[3]} D.A. Goldston, S. Graham, J. Pintz and C.Y. Yildirim: 
Small gaps between primes or almost primes. Trans.\ AMS., {\bf 361}
 (2009), 5285--5330. See also
arXiv: math/0506067 v1. June 2005.
\item{[4]} D.A. Goldston, J. Pintz and C.Y. Yildirim: Primes in tuples.\ I.
Ann.\ Math., (2), {\bf 170} (2009), 819--862.
See also arXiv: math/0508185 v1. August 2005.
\item{[5]} D.A. Goldston and C.Y. Yildirim: Small gaps between primes I.
arXiv: math/\hfill\break0504336 v1. April 2005.
\item{[6]} Yu.V. Linnik: The large sieve. 
C.R. Acad.\ Sci.\ URSS (N.S.),  {\bf 30} (1941), 292--294.
\item{[7]} ---: {\it Dispersion Method in Binary 
Additive Problems\/}. Leningrad Univ.\ Press, Leningrad 1961.
(Russian)
\item{[8]} J. Maynard: Small gaps between primes. arXiv: 1311.4600 v2.
November 2013.
\item{[9]} Y. Motohashi: An induction principle for the
 generalization of  Bombieri's prime number theorem. 
 Proc.\ Japan Acad., {\bf 52} (1976), 273--275.
\item{[10]} ---: {\it Sieve Methods and Prime Number Theory\/}.
Tata IFR\ Lect.\ Math.\ Phy., {\bf 72}, 
Tata IFR--Springer 1983.
\item{[11]} ---: Talk at the AIM workshop 
`{\it Gaps between primes\/}'. November/December 2005. 
http://aimath.org/pastworkshops/primegapsrep.pdf
\item{[12]}---: {\it Analytic Number Theory\/} I. 
Asakura, 2009. (Japanese; English edition is under preparation)
\item{[13]} Y. Motohashi and J. Pintz: A smoothed GPY sieve.
Bull.\ London Math.\ Soc., {\bf 40} (2008), 298--310. See also
arXiv: math/0602599 v1. February 2006;  v2. July 2013.
\item{[14]} J. Pintz: Polignac numbers, conjectures of Erd\"os on
gaps between primes, arithmetic progressions in primes, and the bounded
gap conjecture. arXiv: 1305.6289. May 2013.
\item{[15]} D.H.J. Polymath: New equidistribution estimates of
Zhang type, and bounded gaps between primes. 
arXiv: 1402.0811 v1. February 2014.
\item{[16]} A. R\'enyi: On the representation of an 
even number as the sum of a prime and an almost prime. 
Izv.\ Akad.\ Nauk SSSR Ser.\ Mat.,  {\bf 12} (1948), 
57--78. (Russian)
\item{[17]} A. Selberg:  Lectures on sieves. In {\it Collected Papers\/}.\ II.
Springer, Berlin 1991, pp.\ 65--247.
\item{[18]} T. Tao: http://terrytao.wordpress.com/2013/11/19/
polymath8b-bounded\hfill\break-intervals-with-many-primes-after-maynard/
\item{[19]} R.C. Vaughan: An elementary method in prime 
number theory. Acta Arith.,  {\bf 37} (1980), 111--115.
\item{[20]} A.I. Vinogradov: The density hypothesis for 
Dirichlet $L$-series. Izv.\ Akad.\ Nauk SSSR Ser.\ Mat., 
{\bf 29} (1965), 903--934; Corrigendum. ibid., {\bf 30} (1966), 
719--720. (Russian)
\item{[21]} Y. Zhang: Bounded gaps between primes. Preprint, April 2013.
\medskip
\hfill www.math.cst.nihon-u.ac.jp/$\sim$ymoto/

\bye